\documentclass{amsart}

\usepackage{amsfonts}
\usepackage{amssymb}
\usepackage{amsthm}
\usepackage{eucal}

\theoremstyle{plain} 
\newtheorem{theorem}{Theorem}[section]

\theoremstyle{definition}

\theoremstyle{remark}

\newcommand{\ol}{\overline}

\newcommand{\gen}[1]{\langle#1\rangle}
\newcommand{\parn}[1]{\left(#1\right)}

\newcommand{\vep}{\varepsilon}
\newcommand{\al}{\alpha}
\newcommand{\be}{\beta}
\newcommand{\ga}{\gamma}

\begin{document}
\title{Low dimensional nilpotent $n$-Lie algebras}

\author[M. Eshrati]{mehdi eshrati$^1$}
\author[F. Saeedi]{farshid saeedi$^2$}
\author[H. Darabi]{hamid darabi$^3$}

\date{}

\address{$^{1}$Farhangiyan University, Shahid Beheshti Campus, Mashhad, Iran.}
\address{$^{2}$Department of Mathematics, Mashhad Branch, Islamic Azad University, Mashhad, Iran.}
\address{$^{3}$ Department of Mathematics, Esfarayen University of Technology, Esfarayen, Iran.}

\email{eshrati.me@gmail.com}
\email{saeedi@mshdiau.ac.ir}
\email{darabi@iauesf.ac.ir}

\begin{abstract}
In this paper, nilpotent $n$-Lie algebras of dimension $n+3$ as well as nilpotent $n$-Lie algebras of class $2$ and dimension $n+4$ are classified.
\end{abstract}
\maketitle
\section{Introduction and preliminaries}
The classification of low dimensional Lie algebras is an early problem in the study of Lie algebras. Such classifications can be find in standard references of Lie algebras. A first step to classify $6$-dimensional Lie algebras was done by Umlauf \cite{ku} in 1891. In 1950, Mozorov \cite{vm} obtained the classification of nilpotent Lie algebras of dimension less than $6$ over arbitrary fields and those of dimensions $6$ over a field of characteristic zero. The classification of $6$-dimensional nilpotent Lie algebras is completed by Cicalo \cite{sc.wg.cs} in 2012. The $7$-dimensional nilpotent Lie algebras over the fields of real and complex numbers are classified in \cite{cs}. Also, the $8$-dimensional nilpotent Lie algebras of class $2$ with a $4$-dimensional center, those with a $2$-dimensional center, and those with a $4$-dimensional center over the field of complex numbers are classified in \cite{sx.br}, \cite{br.lz} and \cite{yw.hc.yn}, respectively. We know that a classification of Lie algebras with respect to the dimension of their multiplier is given in \cite{lb}.

In 1985, Fillipov \cite{vf} defined an \textit{$n$-Lie algebra} as an antisymmetric $n$-linear map on a vector space that satisfy the following Jacobi identity:
\[[[x_1,x_2,\ldots,x_n],y_2,\ldots,y_n]=\sum_{i=1}^n[x_1,\ldots,x_{i-1},[x_i,y_2,\ldots,y_n],x_{i+1},\ldots,x_n]\]
for all $x_i,y_i\in L$, $1\leq i\leq n$ and $2\leq j\leq n$. Also, he has classified the $n$-Lie algebras of dimension $n+1$ over an algebraically closed field of characteristic zero. In 2008, Bai \cite{rb.xw.wx.ha} classified those $n$-Lie algebras of dimension $n+1$ whose underlying field has characteristic $2$. Also, in 2011, Bai et. al. \cite{rb.gs.yz} classified the $n$-Lie algebras of dimension $n+2$ over algebraically closed fields of characteristic zero.

Let $A_1,A_2,\ldots,A_n$ be subalgebras of $n$-Lie algebra $A$. Then the subalgebra of $A$ generated by all commutators $[x_1,\ldots,x_n]$, in which $x_i\in A_i$, is denoted by $[A_1,\ldots,A_n]$. The subalgebra $A^2=[A,\ldots,A]$ is called the \textit{derived} $n$-Lie subalgebra of $A$. If $A^2=0$, we call $A$ an abelian $n$-Lie algebra. The \textit{center} of the $n$-Lie algebra $A$ is defined as
\[Z(A)=\{x\in A:[x,A,\ldots,A]=0\}.\]
Assume $Z_0(A)=0$, then $i$th center of $A$ is defined inductively as $Z_i(A)/Z_{i-1}(A)=Z(A/Z_{i-1}(A))$ for all $i\geq1$. In 1987, Kasimov \cite{sk} defines the notion of nilpotency of $n$-Lie algebras as follows:

An $n$-Lie algebra $A$ is called nilpotent if $A^s=0$ for some non-negative integer $s$, in which $A^{i+1}$ is defined inductively as $A^1=A$ and $A^{i+1}=[A^i,A,\ldots,A]$. The $n$-Lie algebra $A$ is nilpotent of class $c$ provided that $A^{c+1}=0$ and $A^c\neq0$. Similar results for the same groups are obtained for nilpotent $n$-Lie algebras in \cite{mw}.

The goal of this paper is to classify $(n+3)$-dimensional nilpotent $n$-Lie algebras as well as $(n+4)$-dimensional nilpotent $n$-Lie algebras of class $2$ over an arbitrary field. This paper is organized as follows: Section 2 presents the preliminary results that will be used in the next sections. In section 3, we give a classification of all $(n+3)$-dimensional nilpotent $n$-Lie algebras over an arbitrary field. Also, in section 4, we shall classify all $(n+4)$-dimensional nilpotent $n$-Lie algebras of class $2$ over an arbitrary field.
\section{Known results}
In this section, we shall present some known results, without proofs, that will be used later. Recall that an $n$-Lie algebra $A$ is called \textit{Special Heisenberg} if $A^2=Z(A)$ and $\dim A^2=1$.
\begin{theorem}[\cite{me.fs.hd}]\label{special Heisenberg}
Every Special Heisenberg $n$-Lie algebra is of dimension $mn+1$ for some positive integer $m$ and it has a presentation as:
\[H(n,m)=\gen{x,x_1,\ldots,x_{mn}:[x_{n(i-1)+1},x_{n(i-1)+2},\ldots,x_{ni}]=x,i=1,\ldots,m}.\]
\end{theorem}
\begin{theorem}[\cite{hd.fs.me}]\label{dimA^2=1}
Let $A$ be a nilpotent $n$-Lie algebra of dimension $d$ satisfying $\dim A^2=1$. Then $A\cong H(n,m)\oplus F(d-mn-1)$ for some $m\geq1$.
\end{theorem}
\begin{theorem}[\cite{hd.fs.me}]\label{A^2=Z(A)}
Let $A$ be a nilpotent $n$-Lie algebra of dimension $d=n+k$, for $3 \leq k\leq n+1$ such that $A^2=Z(A)$ and  $\dim A^2=2$. Then
\[ A\cong \gen{e_1,\ldots,e_{n+k}:[e_{k-1},\ldots,e_{n+k-2}]=e_{n+k},[e_1,\ldots,e_n]=e_{n+k-1}}.\]
This $n$-Lie algebra is denoted by $A_{n,k}$.
\end{theorem}

\begin{theorem}[\cite{hd.fs.me}]\label{d<=n+2}
Let $A$ be a non-abelian nilpotent $n$-Lie algebra of dimension $d\leq n+2$. Then
\[A\cong H(n,1),\ H(n,1)\oplus F(1)\ \text{or}\ A_{n+2,3},\]
in which
\[A_{n+2,3}=\gen{e_1,\ldots,e_{n+2}:[e_1,\ldots,e_n]=e_{n+1},[e_2,\ldots,e_{n+1}]=e_{n+2}}.\]
\end{theorem}
\begin{theorem}[\cite{sc.wg.cs}]
\ 
\begin{itemize}
\item[(1)]
Over a field $F$ of characteristic different from $2$, the list of the isomorphisms type of $6$-dimensional nilpotent Lie algebras is the following: $L_{5,k}\oplus F$ with $k\in\{1,\ldots,9\}$ ; $L_{6,k}$ with $k\in\{10,\ldots,18,20,23,25,\ldots,28\}$; $L_{6,k}(\vep_1)$ with $k\in\{19,21\}$ and $\vep_1\in F^*/(\overset{*}{\mathop \sim})$; $L_{6,k}(\vep_2)$ with $k\in\{22,24\}$ and $\vep_2\in F/(\overset{*}{\mathop \sim})$.

\item[(2)]
Over a field $F$ of characteristic $2$, the isomorphism types of $6$-dimensional nilpotent Lie algebras are:
$L_{5,k}\oplus F$ with $k\in\{1,\ldots,9\}$; $L_{6,k}$ with $k\in\{10,\ldots,18,20,23,25,\ldots,28\}$; $L_{6,k}(\vep_1)$ with $k\in\{19,21\}$ and $\vep_1\in F^*/(\overset{*}{\mathop \sim})$; $L_{6,k}(\vep_2)$ with $k\in\{22,24\}$ and $\vep_2\in F/(\overset{*+}{\mathop \sim})$; $L_{6,k}^{(2)}$ with $k\in\{1,2,5,6\}$; $L_{5,k}^{(2)}(\vep_3)$ with $k\in\{3,4\}$ and $\vep_3\in F^*/(\overset{*+}{\mathop \sim})$; $L_{6,k}^{(2)}(\vep_4)$ with $k\in\{7,8\}$ and $\vep_4\in\{0,\omega\}$.

\end{itemize}
\end{theorem}
\section{Classification of $(n+3)$-dimensional nilpotent $n$-Lie algebras}
By Theorem \ref{d<=n+2}, the number of nilpotent $n$-Lie algebras of dimensions $n$, $n+1$ and $n+2$ are one, two and three up to isomorphism. In this section, we shall classify all $(n+3)$-dimensional nilpotent $n$-Lie algebras over an arbitrary field.

Let $A$ be a nilpotent $n$-Lie algebra of dimension $n+3$ with basis $\{e_1,\ldots,e_{n+3}\}$. If $e_{n+3}$ is a central element of $A$, then $A/\gen{e_{n+3}}$ is a nilpotent $n$-Lie algebra of dimension $n+2$. We discuss on the abelian-ness of $A/\gen{e_{n+3}}$.

If $A/\gen{e_{n+3}}$ is abelian, then brackets in $A$ can be written as
\[[e_1,\ldots,\hat{e}_i,\ldots,\hat{e}_j,\ldots,e_{n+2}]=\theta_{ij}e_{n+3},\quad1\leq i<j\leq n+2.\]
If $\theta_{ij}=0$ for all $1\leq i<j\leq n+2$, then $A$ is an abelian $n$-Lie algebra, which we denote it by $A_{n+3,1}$. If $\theta_{ij}$ are not all equal to zero, then $A$ is non-abelian satisfying $\dim A^2=1$. Hence, by Theorem  \ref{dimA^2=1}, $A\cong H(n,m)\oplus F(n+3-nm-1)$.

In case $n>2$, we must have $m=1$ and hence
\[A=\gen{e_1,\ldots,e_{n+3}:[e_1,\ldots,e_n]=e_{n+3}}\cong H(n,1)\oplus F(2),\]
which we denote it by $A_{n+3,2}$.

Also, in case $n=2$, we have two $n$-Lie algebras for $A$, namely $H(2,2)$ and $H(2,1)\oplus F(2)$.

Now, assume that $A/\gen{e_{n+3}}$ is a non-abelian $n$-Lie algebra. By Theorem \ref{d<=n+2}, we have two possibilities for $A/\gen{e_{n+3}}$:

Case 1: $A/\gen{e_{n+3}}\cong\gen{\ol{e}_1,\ldots,\ol{e}_{n+2}:[\ol{e}_1,\ldots,\ol{e}_n]=\ol{e}_{n+1}}$. Then the brackets in $A$ can be written as
\[\begin{cases}
[e_1,\ldots,e_n]=e_{n+1}+\al e_{n+3},&\\
[e_1,\ldots,\hat{e}_i,\ldots,\hat{e}_j,\ldots,e_n,e_{n+1},e_{n+2}]=\theta_{ij}e_{n+3},&1\leq i<j\leq n,\\
[e_1,\ldots,\hat{e}_i,\ldots,e_n,e_{n+1}]=\be_i e_{n+3},&1\leq i\leq n,\\
[e_1,\ldots,\hat{e}_i,\ldots,e_n,e_{n+2}]=\ga_ie_{n+3},&1\leq i\leq n.
\end{cases}\]

Regarding a suitable change of basis, one can assume that $\al=0$, and the Jacobi identities give us \[\theta_{ij}=0,\quad1\leq i<j\leq n.\]
Hence
\[\begin{cases}
[e_1,\ldots,e_n]=e_{n+1},&\\
[e_1,\ldots,\hat{e}_i,\ldots,e_n,e_{n+1}]=\be_ie_{n+3},&1\leq i\leq n,\\
[e_1,\ldots,\hat{e}_i,\ldots,e_n,e_{n+2}]=\ga_ie_{n+3},&1\leq i\leq n.
\end{cases}\]
The above brackets show that the dimension of the center of $A$ is at most $3$. We discuss on the dimension of the center of $A$.

(i) $\dim Z(A)=1$. Then, we must have $\be_i,\ga_j\neq0$ for some $i$ and $j$. Without loss of generality assume that $\be_1,\ga_1\neq0$. Applying the following transformations
\[e'_1=e_1+\sum_{i=2}^n(-1)^{i-1}\frac{\be_i}{\be_1}e_i,\quad e'_j=e_j,\quad2\leq j\leq n+2,\quad e'_{n+3}=\be_1e_{n+3}\]
we obtain
\[\begin{cases}
[e'_1,\ldots,e'_n]=e'_{n+1},&\\
[e'_2,\ldots,e'_{n+1}]=e'_{n+3},&\\
[e'_2,\ldots,e'_n,e'_{n+2}]=\frac{\ga_1}{\be_1}e'_{n+3},&\\
[e'_1,\ldots,\hat{e}'_i,\ldots,e'_n,e'_{n+2}]=\frac{1}{\be_1}\parn{\ga_i-\frac{\be_i\ga_1}{\be_1}}e'_{n+3},&2\leq i\leq n.
\end{cases}\]
Next, by applying the transformations
\begin{align*}
e''_1&=e'_1+\sum_{i=2}^n(-1)^{i-1}\parn{\frac{\ga_i}{\ga_1}-\frac{\be_i}{\be_1}}e'_i,\\
e''_j&=e'_j,\quad2\leq j\leq n+1,\quad e''_{n+2}=\frac{\be_1}{\ga_1}e'_{n+2},\quad e''_{n+3}=e'_{n+3}
\end{align*}
it yields
\[[e''_1,\ldots,e''_n]=e''_{n+1},\quad[e''_2,\ldots,e''_{n+1}]=e''_{n+3},\quad[e''_2,\ldots,e''_n,e''_{n+2}]=e''_{n+3}.\]
Hence, we conclude that
\[A=\gen{e_1,\ldots,e_{n+3}:[e_1,\ldots,e_n]=e_{n+1},[e_2,\ldots,e_{n+1}]=e_{n+3},[e_2,\ldots,e_n,e_{n+2}]=e_{n+3}},\]
which we denote it by $A_{n+3,3}$.

(ii) $\dim Z(A)=2$. We have two possibilities:

(ii-a) If $Z(A)=\gen{e_{n+2},e_{n+3}}$, then at least one of the $\be_i$ is non-zero while all $\ga_i$ are zero. Without loss of generality assume that $\be_1\neq0$. Then
\[\begin{cases}
[e_1,\ldots,e_n]=e_{n+1},&\\
[e_2,\ldots,e_{n+1}]=\be_1e_{n+3},&\\
[e_1,\ldots,\hat{e}_i,\ldots,e_n,e_{n+1}]=\be_ie_{n+3},&2\leq i\leq n.
\end{cases}\]
Applying the following transformations
\[e'_1=e_1+\sum_{i=2}^n(-1)^{i-1}\frac{\be_i}{\be_1}e_i,\quad e'_j=e_j,\quad2\leq j\leq n+2,\quad e'_{n+3}=\be_1e_{n+3}\]
we obtain
\[\begin{cases}
[e'_1,\ldots,e'_n]=e'_{n+1},&\\
[e'_2,\ldots,e'_{n+1}]=e'_{n+3}.
\end{cases}\]
Hence
\[A=\gen{e_1,\ldots,e_{n+3}:[e_1,\ldots,e_n]=e_{n+1},[e_2,\ldots,e_{n+1}]=e_{n+3}},\]
which is isomorphic to $A_{n+2,3}\oplus F(1)$ and it is denoted by $A_{n+3,4}$.

(ii-b) If $Z(A)=\gen{e_{n+1},e_{n+3}}$, then $A^2=Z(A)$ and Theorem \ref{A^2=Z(A)} yields
\[A=\gen{e_1,\ldots,e_{n+3}:[e_1,\ldots,e_n]=e_{n+1},[e_2,\ldots,e_n,e_{n+2}]=e_{n+3}},\]
which is denoted by $A_{n+3,5}$.

(iii) $\dim Z(A)=3$. Then $Z(A)=\gen{e_{n+1},e_{n+2},e_{n+3}}$ and so that $\be_i=\ga_i=0$ for all $1\leq i\leq n$. The only non-zero bracket is $[e_1,\ldots,e_n]=e_{n+1}$, which gives rise to the algebra $H(n,1)\oplus F(2)=A_{n+3,2}$.

Case 2: $A/\gen{e_{n+3}}\cong\gen{\ol{e}_1,\ldots,\ol{e}_{n+2}:[\ol{e}_1,\ldots,\ol{e}_n]=\ol{e}_{n+1},[\ol{e}_2,\ldots,\ol{e}_{n+1}]=\ol{e}_{n+2}}$. Then the brackets in $A$ are as follows:
\[\begin{cases}
[e_1,\ldots,e_n]=e_{n+1}+\al e_{n+3},&\\
[e_2,\ldots,e_{n+1}]=e_{n+2}+\be  e_{n+3},&\\
[e_1,\ldots,\hat{e}_i,\ldots,\hat{e}_j,\ldots,e_n,e_{n+1},e_{n+2}]=\theta_{ij}e_{n+3},&1\leq i<j\leq n,\\
[e_2,\ldots,e_n,e_{n+2}]=\gamma e_{n+3},&\\
[e_1,\ldots,\hat{e}_i,\ldots,e_n,e_{n+2}]=\al_ie_{n+3},&2\leq i\leq n,\\
[e_1,\ldots,\hat{e}_i,\ldots,e_n,e_{n+1}]=\be_ie_{n+3},&2\leq i\leq n.
\end{cases}\]
With a suitable change of basis, one can assume that $\al=\be=0$. Moreover, from the Jacobi identities, it follows that
\[\theta_{ij}=0,\quad1\leq i<j\leq n,\quad\al_i=0,\quad2\leq i\leq n,\]
so
\[\begin{cases}
[e_1,\ldots,e_n]=e_{n+1},&\\
[e_2,\ldots,e_{n+1}]=e_{n+2},&\\
[e_2,\ldots,e_n,e_{n+2}]=\gamma e_{n+3},&\\
[e_1,\ldots,\hat{e}_i,\ldots,e_n,e_{n+1}]=\be_ie_{n+3},&2\leq i\leq n.
\end{cases}\]
The above relations show that the dimension of $Z(A)$ is at most $2$. We have two possibilities:

(i) $\dim Z(A)=1$. Then $\ga\neq0$ and $e'_{n+3}=\ga e_{n+3},\quad e'_j=e_j,\quad 1\leq j\leq n+2$ yields
\[\begin{cases}
[e'_1,\ldots,e'_n]=e'_{n+1},&\\
[e'_2,\ldots,e'_{n+1}]=e'_{n+2},&\\
[e'_2,\ldots,e'_n,e'_{n+2}]=e'_{n+3},&\\
[e'_1,\ldots,\hat{e'}_i,\ldots,e'_n,e'_{n+1}]=\frac{\be_i}{\ga}e'_{n+3},&2\leq i\leq n.
\end{cases}\]
If $\be_i=0$ for all $2\leq i\leq n$, then
\[A=\gen{e_1,\ldots,e_{n+3}:[e_1,\ldots,e_n]=e_{n+1},[e_2,\ldots,e_{n+1}]=e_{n+2},[e_2,\ldots,e_n,e_{n+2}]=e_{n+3}},\]
which is denoted by $A_{n+3,6}$. On the other hand, if $\be_i\neq0$ for some $2\leq i\leq n$, say $\be_2\neq0$, then, by applying the following transformations,
\begin{align*}
e'_1&=\frac{\be_2}{\ga}e_1,\quad e'_2=\parn{\frac{\be_2}{\ga}}^2 e_2+\sum_{i=3}^n(-1)^i\frac{\be_2\be_i}{\ga}e_i,\quad e'_i=e_i,\quad3\leq i\leq n-1,\\
e'_n&=\parn{\frac{\ga}{\be_2}}^2e_n,e'_i=\frac{\be_2}{\ga}e_i,\quad n+1\leq i\leq n+2,\quad e'_{n+3}=\frac{\be_2}{\ga}e_{n+3},
\end{align*}
we observe that
\begin{multline*}
A=\langle e_1,\ldots,e_{n+3}:[e_1,\ldots,e_n]=e_{n+1},[e_2,\ldots,e_{n+1}]=e_{n+2},\\
[e_2,\ldots,e_n,e_{n+2}]=[e_1,e_3,\ldots,e_{n+1}]=e_{n+3}\rangle,
\end{multline*}
which we denote it by $A_{n+3,7}$.

(ii) $\dim Z(A)=2$. Then $\ga=0$ and the brackets in $A$ reduce to
\[\begin{cases}
[e_1,\ldots,e_n]=e_{n+1},&\\
[e_2,\ldots,e_{n+1}]=e_{n+2},&\\
[e_1,\ldots,\hat{e}_i,\ldots,e_n,e_{n+1}]=\be_ie_{n+3},&2\leq i\leq n.
\end{cases}\]
If $\be_i=0$ for all $2\leq i\leq n$, then
\[A=\gen{e_1,\ldots,e_{n+3}:[e_1,\ldots,e_n]=e_{n+1},[e_2,\ldots,e_{n+1}]=e_{n+2}}.\]
One can easily see that this algebra is isomorphic to $A_{n+3,4}$, while if $\be_i\neq0$ for some $2\leq i\leq n$, say $\be_2\neq0$, then the following transformations
\[e'_1=e_1,\quad e'_2=e_2+\sum_{i=3}^n(-1)^i\frac{\be_i}{\be_2}e_i,\quad e'_j=e_j,\quad3\leq j\leq n+2,\quad e'_{n+3}=\be_2e_{n+3}\]
show  that
\begin{multline*}
A=\langle e_1,\ldots,e_{n+3}:[e_1,\ldots,e_n]=e_{n+1},\\
[e_2,\ldots,e_{n+1}]=e_{n+2},[e_1,e_3,\ldots,e_{n+1}]=e_{n+3}\rangle.
\end{multline*}
This algebra is denoted by $A_{n+3,8}$.

The above results summary as:
\begin{theorem}\label{d=n+3}
The only nilpotent $n$-Lie algebras of dimension $n+3$ with $n>2$ are $A_{n+3,i}$ with $i\in\{1,\ldots,8\}$. For $n=2$ we have instead the algebra $H(2,2)$.
\end{theorem}
\section{Classification of $(n+4)$-dimensional nilpotent $n$-Lie algebras of class $2$}
Nilpotent $n$-Lie algebras of class $2$ appear in some problems of geometry like commutative Riemannian manifolds. Also, classifying nilpotent Lie algebras of class $2$ has been an important problem in Lie algebras. In \cite{br.lz} nilpotent Lie algebras of class $2$ and dimension $8$ with $2$-dimensional center over the field of complex numbers are classified. In this section, we intend to classify $(n+4)$-dimensional nilpotent $n$-Lie algebras of class $2$. Regarding Theorem \ref{d=n+3}, the following result obtains immediately.
\begin{theorem}\label{d=n+3 of class 2}
The only $(n+3)$-dimensional nilpotent $n$-Lie algebras of class $2$ are
\begin{itemize}
\item[(1)]$A_{n+3,2}=\gen{e_1,\ldots,e_{n+3}:[e_1,\ldots,e_n]=e_{n+3}}\cong H(n,1)\oplus F(2)$;
\item[(2)]$A_{n+3,5}=\gen{e_1,\ldots,e_{n+3}:[e_1,\ldots,e_n]=e_{n+1},[e_2,\ldots,e_n,e_{n+2}]=e_{n+3}}$;
\item[(3)]$H(2,2)$.
\end{itemize}
\end{theorem}

Suppose $A$ is an $(n+4)$-dimensional nilpotent $n$-Lie algebra of class $2$ with basis $\{e_1,\ldots,e_{n+4}\}$. If $\dim A^2=1$, Theorem  \ref{dimA^2=1} yields $A\cong H(n,m)\oplus F(n+4-nm-1)$. Hence $A$ is isomorphic to one of the following algebras
\[H(n,1)\oplus F(3),\quad H(2,2)\oplus F(1)\ \text{or}\ H(3,2).\]
Now, we assume that $\dim A^2\geq2$ and $\gen{e_{n+3},e_{n+4}}\subseteq A^2$. Then $A/\gen{e_{n+4}}$ is an $(n+3)$-dimensional nilpotent $n$-Lie algebra of class $2$ so that, by Theorem \ref{d=n+3 of class 2}, $A/\gen{e_{n+4}}$ is of the following forms:

Case 1: $A/\gen{e_{n+4}}\cong\gen{\ol{e}_1,\ldots,\ol{e}_{n+3}:[\ol{e}_1,\ldots,\ol{e}_n]=\ol{e}_{n+3}}$. Then the brackets in $A$ can be written as:
\[\begin{cases}
[e_1,\ldots,e_n]=e_{n+3}+\al e_{n+4},&\\
[e_1,\ldots,\hat{e}_i,\ldots,\hat{e}_j,\ldots,e_n,e_{n+1},e_{n+2}]=\theta_{ij}e_{n+4},&1\leq i<j\leq n,\\
[e_1,\ldots,\hat{e}_i,\ldots,e_n,e_{n+1}]=\be_ie_{n+4},&1\leq i\leq n,\\
[e_1,\ldots,\hat{e}_i,\ldots,e_n,e_{n+2}]=\ga_ie_{n+4},&1\leq i\leq n.
\end{cases}\]
By a suitable change of basis, we may assume that $\al=0$, hence
\[\begin{cases}
[e_1,\ldots,e_n]=e_{n+3},&\\
[e_1,\ldots,\hat{e}_i,\ldots,\hat{e}_j,\ldots,e_n,e_{n+1},e_{n+2}]=\theta_{ij}e_{n+4},&1\leq i<j\leq n,\\
[e_1,\ldots,\hat{e}_i,\ldots,e_n,e_{n+1}]=\be_ie_{n+4},&1\leq i\leq n,\\
[e_1,\ldots,\hat{e}_i,\ldots,e_n,e_{n+2}]=\ga_ie_{n+4},&1\leq i\leq n.
\end{cases}\]
Since $\dim A^2\geq2$, it is obvious that $\dim Z(A)\leq3$. First, assume that $\dim Z(A)=3$. Then, without loss of generality, we may assume that $Z(A)=\gen{e_{n+2},e_{n+3},e_{n+4}}$. Hence $\ga_i$ and $\theta_{ij}$ are all zero while $\be_i\neq0$ for some $i$. We may assume that $\be_1\neq0$. Then
\[\begin{cases}
[e_1,\ldots,e_n]=e_{n+3},&\\
[e_2,\ldots,,e_n,e_{n+1}]=\be_1e_{n+4},&\\
[e_1,\ldots,\hat{e}_i,\ldots,e_n,e_{n+1}]=\be_ie_{n+4},&2\leq i\leq n.
\end{cases}\]
Using the following transformations
\[e'_1=e_1+\sum_{i=1}^n(-1)^{i-1}\frac{\be_i}{\be_1}e_i,\quad e'_j=e_j,\quad2\leq j\leq n+3,\quad e'_{n+4}=\be_1e_{n+4},\]
it follows that
\[A=\gen{e_1,\ldots,e_{n+4}:[e_1,\ldots,e_n]=e_{n+3},[e_2,\ldots,e_{n+1}]=e_{n+4}}.\]
This algebra is denoted by $A_{n+4,1}$.

Next, assume that $\dim Z(A)=2$. Then $Z(A)=A^2$. If $n\geq3$, then, by Theorem \ref{A^2=Z(A)}, we have
\[A=\gen{e_1,\ldots,e_{n+4}:[e_1,\ldots,e_n]=e_{n+3},[e_3,\ldots,e_{n+2}]=e_{n+4}},\]
which is denoted by $A_{n+3,2}$. From \cite{sc.wg.cs} in case $n=2$, the only $6$-dimensional Lie algebra satisfying $A^2=Z(A)$ and $\dim Z(A)=2$ is
\[A\cong L_{6,22}(\vep)=\gen{e_1,\ldots,e_6:[e_1,e_2]=e_5,[e_1,e_3]=e_6,[e_2,e_4]=\vep e_6,[e_3,e_4]=e_5}.\]
Note that this algebra does not satisfy $A/\gen{e_6}\cong H(2,1)\oplus F(2)$.

Case 2: $A/\gen{e_{n+4}}\cong\gen{\ol{e}_1,\ldots,\ol{e}_{n+3}:[\ol{e}_1,\ldots,\ol{e}_n]=\ol{e}_{n+1},[\ol{e}_2,\ldots,\ol{e}_n,\ol{e}_{n+2}]=\ol{e}_{n+3}}$. Then the brackets in $A$ can be written as
\[\begin{cases}
[e_1,\ldots,e_n]=e_{n+1}+\al e_{n+4},&\\
[e_2,\ldots,e_n,e_{n+2}]=e_{n+3}+\be e_{n+4},&\\
[e_1,\ldots,\hat{e}_i,\ldots,\hat{e}_j,\ldots,e_n,e_{n+1},e_{n+2}]=\theta_{ij}e_{n+4},&1\leq i<j\leq n,\\
[e_1,\ldots,\hat{e}_i,\ldots,e_n,e_{n+1}]=\be_ie_{n+4},&1\leq i\leq n,\\
[e_1,\ldots,\hat{e}_i,\ldots,e_n,e_{n+2}]=\ga_ie_{n+4},&2\leq i\leq n.
\end{cases}\]
With a suitable change of basis, one can assume that $\al=\be=0$. A simple verification shows $Z(A)=\gen{e_{n+1},e_{n+3},e_{n+4}}$ such that
\[\begin{cases}
[e_1,\ldots,e_n]=e_{n+1},&\\
[e_2,\ldots,e_n,e_{n+2}]=e_{n+3},&\\
[e_1,\ldots,\hat{e}_i,\ldots,e_n,e_{n+2}]=\ga_ie_{n+4},&2\leq i\leq n.
\end{cases}\]
If $\ga_i=0$ for all $2\leq i\leq n$, then
\[A=\gen{e_1,\ldots,e_{n+4}:[e_1,\ldots,e_n]=e_{n+1},[e_2,\ldots,e_n,e_{n+2}]=e_{n+3}}.\]
One can easily see that this algebra is isomorphic to $A_{n+4,1}$. On the other hand, if $\ga_i\neq0$ for some $2\leq i\leq n$, say $\ga_2\neq0$, then we may apply the following transformations
\[e'_2=e_2+\sum_{j=3}^n(-1)^{j-1}\frac{\ga_j}{\ga_2}e_j,\quad e'_i=e_i,\quad i=1,3,\ldots,n+3,\quad e'_{n+4}=\ga_2e_{n+4}.\]
Then
\begin{multline*}
A=\langle e_1,\ldots,e_{n+4}:[e_1,\ldots,e_n]=e_{n+1},\\
[e_2,\ldots,e_n,e_{n+2}]=e_{n+3},[e_1,e_3,\ldots,e_n,e_{n+2}]=e_{n+4}\rangle,
\end{multline*}
which we denote it by $A_{n+4,3}$.

Case 3: $A/\gen{e_{n+4}}\cong H(2,2)$. From \cite{sc.wg.cs}, we observe that 
\[A\cong L_{6,22}(\vep)=\gen{e_1,\ldots,e_6:[e_1,e_2]=e_5,[e_1,e_3]=e_6,[e_2,e_4]=\vep e_6,[e_3,e_4]=e_5}.\]
\begin{theorem}
The only $(n+4)$-dimensional nilpotent $n$-Lie algebras of class $2$ are
\[H(n,1)\oplus F(3),\ H(2,2)\oplus F(1),\ H(3,2),\ A_{n+4,1},\ A_{n+4,2},\ A_{n+4,3}\ \text{and}\ L_{6,22}(\vep).\]
\end{theorem}

In table I, we have illustrated all $n$-Lie algebras obtained in this paper.
\begin{center}
Table I
\begin{tabular}{|c|c|}
\hline
Nilpotent $n$-Lie algebra&Non-zero multiplications\\
\hline
$A_{n+3,1}$&--\\
\hline
$A_{n+3,2}$&$[e_1,\ldots,e_n]=e_{n+3}$\\
\hline
$A_{n+3,3}$&$\begin{array}{c}[e_1,\ldots,e_n]=e_{n+1},[e_2,\ldots,e_{n+1}]=e_{n+3},\\{[e_2,\ldots,e_n,e_{n+2}]}=e_{n+3}\end{array}$\\
\hline
$A_{n+3,4}$&$[e_1,\ldots,e_n]=e_{n+1},[e_2,\ldots,e_{n+1}]=e_{n+3}$\\
\hline
$A_{n+3,5}$&$[e_1,\ldots,e_n]=e_{n+1},[e_2,\ldots,e_n,e_{n+2}]=e_{n+3}$\\
\hline
$A_{n+3,6}$&$\begin{array}{c}[e_1,\ldots,e_n]=e_{n+1},[e_2,\ldots,e_{n+1}]=e_{n+2},\\{[e_2,\ldots,e_n,e_{n+2}]}=e_{n+3}\end{array}$\\
\hline
$A_{n+3,7}$&$\begin{array}{c}[e_1,\ldots,e_n]=e_{n+1},[e_2,\ldots,e_{n+1}]=e_{n+2},\\{[e_2,\ldots,e_n,e_{n+2}]}=[e_1,e_3,\ldots,e_{n+1}]=e_{n+3}\end{array}$\\
\hline
$A_{n+3,8}$&$\begin{array}{c}[e_1,\ldots,e_n]=e_{n+1},[e_2,\ldots,e_{n+1}]=e_{n+2},\\{[e_1,e_3,\ldots,e_{n+1}]}=e_{n+3}\end{array}$\\
\hline
$A_{n+4,1}$&$[e_1,\ldots,e_n]=e_{n+3},[e_2,\ldots,e_{n+1}]=e_{n+4}$\\
\hline
$A_{n+4,2}$&$[e_1,\ldots,e_n]=e_{n+3},[e_3,\ldots,e_{n+2}]=e_{n+4}$ ($n\geq3$)\\
\hline
$A_{n+4,3}$&$\begin{array}{c}[e_1,\ldots,e_n]=e_{n+1},[e_2,\ldots,e_n,e_{n+2}]=e_{n+3},\\{[e_1,e_3,\ldots,e_n,e_{n+2}]}=e_{n+4}\end{array}$\\
\hline
\end{tabular}
\end{center}

\end{document}